\documentclass[5p,preprint]{elsarticle}

\usepackage{amsmath,amsfonts,amsthm,dsfont,amssymb}
\usepackage{mathdots}
\usepackage{mdwlist}
\usepackage{mathrsfs}
\usepackage{comment}
\usepackage{setspace}
\usepackage{verbatim}
\usepackage{cancel}

\usepackage{float}
\usepackage{longtable}
\usepackage{supertabular}
\usepackage{subfig}
\usepackage{bm}
\usepackage{overpic}
\usepackage{enumerate}
\usepackage{siunitx}
\usepackage[europeanresistors]{circuitikz}

\usepackage{tikz,pgfplots}
\usetikzlibrary{fit}
\usetikzlibrary{calc}
\usetikzlibrary{arrows}
\pgfdeclarelayer{background}
\colorlet{inbox}{lightgray!20}
\colorlet{outbox}{lightgray!50}

\usetikzlibrary{patterns,math}
\usetikzlibrary{arrows}
\usetikzlibrary{shapes}
\usetikzlibrary{decorations.pathmorphing}
\usetikzlibrary{decorations.pathreplacing}
\usetikzlibrary{decorations.shapes}
\usetikzlibrary{decorations.text}
\usetikzlibrary{patterns}
\usetikzlibrary{backgrounds}
\usepackage{tkz-euclide}

  \pgfdeclarelayer{background1}
  \pgfdeclarelayer{background2}
  \pgfdeclarelayer{background3}
  \pgfdeclarelayer{background4}
  \pgfdeclarelayer{foreground}
  \pgfsetlayers{background4,background3,background2,background1,foreground,background,main}
  \tikzset{errorstyle/.style={thick,red,solid}}
  \tikzset{yrefstyle/.style={thick,black,dashed}}
  \tikzset{safetystyle/.style={thick,blue,dotted}}
  \tikzset{funnelstyle/.style={thick,blue,densely dotted}}
  \tikzset{funnelbackground/.style={black!20,opacity=0.5}}
  \tikzset{funneldstyle/.style={thin,blue,dashed}}
  \tikzset{qstyle/.style={green!50!black,ultra thick}}
  \tikzset{qhelpstyle/.style={green!50!black,very thin}}
  \tikzset{funnelfillstyle/.style={blue!20!white,opacity=0.8}}
  \tikzset{safetyfillstyle/.style={blue,opacity=0.1}}
  \tikzset{funnelthinfillstyle/.style={blue!5!white,opacity=0.8}}
  \tikzset{safetythinfillstyle/.style={blue!50!white,opacity=0.1}}

\usepackage[colorlinks]{hyperref}

\newtheorem{thm}{Theorem}[section]
 \newtheorem{cor}[thm]{Corollary}
 
 \newtheorem{prop}[thm]{Proposition}
 \theoremstyle{definition}
 \newtheorem{defn}[thm]{Definition}

 \numberwithin{equation}{section}

\newcommand{\cF}{\mathcal{F}}

\newcommand{\half}{\tfrac{1}{2}}

\newcommand{\quarter}{\tfrac{1}{4}}

\newcommand{\N}{\mathbb{N}}

\newcommand{\R}{\mathbb{R}}

\newcommand{\dd}{\textrm d}

\renewcommand{\imath}{\mathrm{i}}

\numberwithin{equation}{section}

{\left(\begin{smallmatrix}}
{\end{smallmatrix}\right)}

{\left[\begin{smallmatrix}}
{\end{smallmatrix}\right]}

\begin{document}

\begin{frontmatter}

\title{Some Remarks on Positive/Negative Feedback}

\tnotetext[funding]{T. Berger acknowledges funding by the Deutsche Forschungsgemeinschaft (DFG, German Research Foundation) -- Project-ID 544702565.}

\author[1]{Thomas Berger}\ead{thomas.berger@mathematik.uni-halle.de}
\author[2]{Achim Ilchmann}\ead{achim.ilchmann@tu-ilmenau.de}
\author[3]{Eugene~P.~Ryan}\ead{masepr@bath.ac.uk}
\address[1]{Institut f\"ur Mathematik, Martin-Luther-Universit\"at Halle-Wittenberg, Theodor-Lieser-Stra{\ss}e 5,  06120 Halle, Germany}
\address[2]{Institut f\"ur Mathematik, Technische Universit\"{a}t Ilmenau, Weimarer Stra{\ss}e 25, 98693~Ilmenau, Germany}
\address[3]{Department of Mathematical Sciences, University of Bath, Bath BA2 7AY, United Kingdom}
\begin{abstract}
In the context of linear control systems, a commonly-held intuition is that
negative and positive feedback cannot both be stability enhancing.  The canonical linear prototype is the
scalar system $\dot x=u$ which, under negative linear feedback $u=-kx$ ($k >0$) is exponentially stable
for all $k >0 $, whereas the lack of exponential instability of the (marginally stable) uncontrolled system is amplified by positive feedback $u=kx$ ($k >0)$.  
By contrast,
for nonlinear systems it is shown, by example, that this intuitive dichotomy may fail to hold.
\end{abstract}

\end{frontmatter}

\section{Nonlinear scalar system}
Consider a scalar system, with state $x$, control $u$, and (unknown) bounded perturbation $p$, of the form
\begin{equation}\label{sys-gen}
\dot x(t)=G(p(t),x(t),u(t)),\quad x(0)=x^0\in\R,
\end{equation}
where $G\colon \R^3\to\R$ and $p\colon\R_{\ge 0}\to\R$  are continuous.  In the unperturbed linear case,
\begin{equation}\label{lin}
\dot x(t)=ax(t)+bu(t),\quad x(0)=x^0, \quad a,b\in\R,~b\ne 0,
\end{equation}
the system implicitly assigns (via the parameter $b\ne 0$) a {\em control direction} to the input, 
namely, $\text{\rm{sgn}}(b)$.  A linear feedback $u(t)=kx(t)$, with $k\ne 0$, is deemed to be {\em positive} 
(respectively, {\em negative}) if $\text{\rm{sgn}}(k)=\text{\rm{sgn}}(b)$ (respectively, $\text{\rm{sgn}}(k)=-\text{\rm{sgn}}(b)$).  In words, 
the feedback is positive if the polarity of the feedback gain coincides with the control direction (and so $bu(t)=+|bk|x(t)$) and is negative 
if the polarity of the gain is opposite to the control direction (and so $u(t)=-|bk|x(t)$).  Whilst the concept of positive/negative feedback is 
unambiguous in the linear case, it is less clear in a nonlinear setting.  We seek a nonlinear generalization.  There is a surfeit of candidates for 
this role.  We focus on one such candidate based on the observation that, for every compact $K\subset \R$ and $v^*\in (0,1)$, the function 
\[
\chi\colon\R\to\R,\quad s\mapsto \min \{v(az+bsv)|~z\in K, ~v^* \le |v|\le 1\}
\]
is such that 
\begin{multline*}
\sup_{s\ge 0} \chi(\eta s)=+\infty \iff
\\
\eta ~\text{coincides with the control direction}~\text{\rm{sgn}}(b).
\end{multline*}
Extrapolating this observation to the nonlinear context of \eqref{sys-gen}, we adopt the following definition of control direction.  
\begin{defn}\label{cd}
$\eta\in\{-1,+1\}$ is a {\em control direction} for \eqref{sys-gen} if, for every compact $K\subset \R\times\R$, there exists $v^*\in (0,1)$ such that the 
function
\begin{multline}\label{chi}
\chi\colon \R\to\R,
\\
s\mapsto\min \{vG(\rho,\xi,sv)|~(\rho,\xi)\in K,~v^*\le |v|\le 1\}
\end{multline}
has the property $\sup_{s\ge 0}\chi (\eta s)=+\infty$.
\end{defn}
There are three possibilities: system \eqref{sys-gen} may fail to have a control direction 
(as is the case when $G$ is bounded); \eqref{sys-gen} has unique control direction (as is so when $G$ is linear);  
\eqref{sys-gen} has non-unique control direction (in Section 3 below, we construct one such example for which each of $\eta=-1$ and $\eta=+1$ qualifies  
as a control direction).

In essence, Definition \ref{cd} attempts to capture the following intuitive idea (loosely stated).  The existence of a control direction is a system requirement that large 
input values $v$ should dominate in the sense that, with the pair $(\rho,\xi)$ of its first two arguments compactly constrained, the function $G$ is such that 
$|G(\rho,\xi,v)|$ can attain arbitrarily large values by taking $|v|$ sufficiently large.  For example, if $G\colon (\rho,\xi,v)\mapsto f(\rho,\xi)g(v)$, then surjectivity of 
(continuous) $g$ is a necessary condition for the existence of a control direction.

For later use, we record the following fact (a particular consequence of a general result 
in variational analysis \cite[Theorem 1.17]{0}).

\begin{prop}\label{prop-cont}
Let $G\colon\R^3\to\R$ be continuous, $K\subset \R\times\R$ compact and $v^*\in (0,1)$.  The function $\chi$, given by \eqref{chi}, is continuous.
\end{prop}

Assuming that \eqref{sys-gen} has a control direction $\eta$, the objective is to determine a continuous  function 
$h\colon\R_{\ge 0}\times\R \to\R$ such that, for every $x^0\in\R$ and bounded, continuous  $p$,  
application  of the feedback $u(t)= -\eta h(t,x(t))$ yields a feedback-controlled initial-value problem
which has a solution $x(\cdot)$, every solution has a global extension to a solution on  $\R_{\ge 0}$, every global solution is such that $x(t)\to 0$ as $t\to\infty$, 
and the control input $u\colon t\mapsto -\eta h(t,x(t))$ is bounded.  

If $\eta=+1$ (respectively, $\eta =-1$), then the feedback $u(t)=-h(t,x(t))$ (respectively, $u(t)=h(t,x(t))$) is deemed to be {\em negative} (respectively, {\em positive}).  
If negative feedback ensures 
benign behaviour in the form of convergence of the state to zero, then intuition might suggest that positive feedback 
causes malign behaviour in the form of non-convergence of the state to zero.  There are many circumstances wherein such a benign/malign dichotomy is 
correct, as discussed in, {\em inter alia}, \cite{1}. 

In Section 3, we show, by a simple example, that there are also circumstances wherein the dichotomy fails to hold: 
a feedback and its reversed-polarity  counterpart can both be benign. 
\section{Feedback system}
Let $\varphi\colon \R_{\ge 0}\to\R_{\ge 0}$ be a continuously differentiable bijection with the property that, for some $c_\varphi >0$,
$\dot\varphi (t)\le c_\varphi\big(1+\varphi(t)\big)$ for all $t\ge 0$ (the simplest example being the identity function with $c_\varphi =1$).
Define the set
\begin{equation}\label{F}\cF:=\{(t,\xi )\in\R_{\ge 0}\times\R\mid ~\varphi(t)|\xi| < 1\},
\end{equation}
the continuous function
\begin{equation}\label{h}
\begin{aligned}
&h\colon \cF\to\R,~~(t,\xi)\mapsto \alpha (\varphi(t)|\xi|)\varphi(t)\xi,\\
&\text{with}~\alpha\colon [0,1)\to [1,\infty),~s\mapsto \frac{1}{1-s},
\end{aligned}
\end{equation}
and the feedback-controlled initial-value problem, on the domain $\cF$,
\begin{equation}\label{fghivp}
\dot x(t)=G(p(t),x(t),-\eta h(t,x(t)),~~x(0)=x^0,
\end{equation}
with continuous $p:\R_{\ge 0}\to \R$.
By a solution we mean a continuous function $x$ on a left-closed interval $I$ (with left endpoint $0$) and $\text{graph}(x)\subset\cF$, satisfying \eqref{fghivp}.
A solution is {\em maximal} if it has no proper right extension that is also a solution.  By the standard theory of ordinary differential equations, for each $x^0\in\R$,
there exists a solution, and every solution has a maximal extension; moreover, the closure of the graph of a maximal solution is not a compact subset of $\cF$.
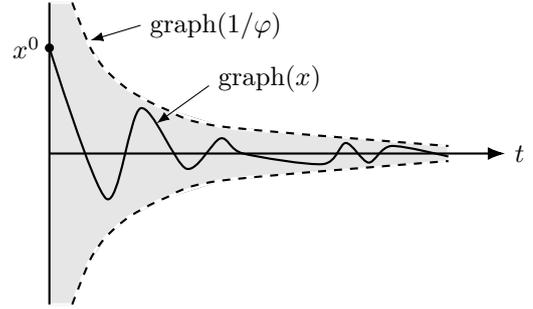
\begin{figure}[h]
 \centering
        \begin{tikzpicture}[x=3cm, y=2cm]

            \definecolor{lightlightgray}{gray}{.9}
            \definecolor{lightergray}{gray}{0.8}
           \filldraw[fill=lightlightgray] {(0,1)--(.1,1) --(0.18,.7)--(.25,.55)--(.35,.42)-- (.5,.3)-- (.75,.18) -- (1,.145) -- (1.75,.05)--(1.75,-.05)--(1,-.145)--(.75,-.18)--(.5,-.3)--(.35,-.42)--(.25,-.55)--(.18,-.7)--(.1,-1)--(0,-1)};
    \draw[thick,white](0,1)--(0.1,1); \draw[thick,white] (1,0.145) -- (1.75,0.05);\draw[thick,white] (1,-0.145) -- (1.75,-0.05);
     \draw[thick,white](0,-1)--(0.1,-1);\draw[thick,white](1.75,.15)--(1.75,-.15);
     \draw[very thick,white]  plot [smooth] coordinates {(.1,1)(.18,.7) (0.25,.55) (.35,.42)(0.5,.3)(.75,.18)(1,.145)};
      \draw[very thick,white]  plot [smooth] coordinates {(.1,-1) (.18,-.7)(0.25,-.55) (.35,-.42)(0.5,-.3)(.75,-.18)(1,-.145)};
            \tikzstyle{dot}=[circle, fill, inner sep=1.2pt]
            \draw[thick,->,arrows=-Latex] (0,0) -- (2,0);\draw[thick] (0,-1) -- (0,1);
             \draw[thick,dashed]  plot [smooth] coordinates {(.1,1)(.18,.7) (0.25,.55)(.35,.42) (0.5,.3)(.75,.18)(1.75,.05)};
              \draw[thick,dashed]  plot [smooth] coordinates {(.1,-1)(.18,-.7) (0.25,-.55)(.35,-.42) (0.5,-.3)(.75,-.18)(1.75,-.05)};
              \draw[thick]  plot [smooth] coordinates{(0,.7)(.25,-.3)(.4,.3)(.6,-.1)(.75,.1)(.85,0)(1.2,-.07)(1.3,.07)(1.4,-.06)(1.5,.05)(1.75,-.02)};
               \node[right] at (2,0){$t$};
               \draw (0,.7) node[dot] {};\node[left] at (0,.7){$x^0$};\node[right] at (.7,.5){$\text{graph}(x)$};\draw[thin,->,arrows=-Latex] (.7,.45) -- (.47,.2);
               \node[right] at (.4,.85){$\text{graph}(1/\varphi)$};\draw[thin,->,arrows=-Latex] (.4,.85) -- (.18,.75);
                             \end{tikzpicture}
    \caption{Domain $\mathcal{F}$.} \label{Fig:scalar-funnel}
\end{figure}
Loosely speaking, the r\^ole of the feedback is to maintain the system's evolution away from the boundary of the domain $\cF$ (wherein singularity resides),
in which case the evolution
continues indefinitely to the right, with transient and asymptotic behaviour determined by the choice of $\varphi$. 

The above feedback is based on a wider body of work in the area of ``funnel control'' (a terminology consistent with Figure \ref{Fig:scalar-funnel}).  
The interested reader is referred to \cite{2} for a review of this area.

\begin{prop}\label{prop}
Assume that $G\colon\R^3\to\R$ is continuous and such that \eqref{sys-gen} has a control direction $\eta\in\{-1,+1\}$.  Let $h$ be as in \eqref{h}.   
For each $x^0\in\R$ and every bounded, continuous $p:\R_{\ge 0}\to \R$, \eqref{fghivp} has a solution and every solution 
has a global extension.  Every global solution
$x\colon\R_{\ge 0}\to \R$  is such that $x(t)\to 0$ as $t\to\infty$, and the control function
$u\colon t\mapsto -\eta h(t,x(t))$ is bounded.
\end{prop}
\begin{proof}
Let $x^0\in \R$ and bounded, continuous $p:\R_{\ge 0}\to \R$ be arbitrary.
Let compact $P\subset\R$ be such that $p(t)\in P$ for all $t\ge 0$.
As already noted, by the standard theory of ordinary differential equations, we know that  \eqref{fghivp} has a solution
and every solution can be maximally extended.   Let $x\colon [0,\omega)\to\R$ be a maximal solution.
It suffices to show existence of ~$\varepsilon > 0$ such that $\varphi(t)|x(t)| < 1-\varepsilon$
for all~$t\in [0,\omega)$, in which case it immediately follows that $\omega =\infty$ (as, otherwise, $\text{graph}(x)$ has compact closure in~$\cF$
which is impossible), $x(t)\to 0$ as $t\to\infty$ (since the bijection $\varphi$ is unbounded), and $|u(t)|=|h(t,x(t))|\le (1-\varepsilon)\alpha (1-\varepsilon)<\infty$ for all $t\ge 0$.
The task of establishing the existence of such ~$\varepsilon >0$ is, in turn, equivalent to proving
boundedness of the function $k\colon t\mapsto \alpha(\varphi(t)|x(t)|)$.   This we proceed to show.

Choose $\tau\in (0,\omega)$ arbitrarily, fix $\mu >1/\varphi(\tau)$ such that  $\mu\ge \max_{t\in [0,\tau]}|x(t)|$, and write
$\varphi_\mu\colon t\mapsto \max\{\mu^{-1},\varphi (t)\}$.  Since $\text{graph}(x)\in\cF$ and $\varphi$ (bijective) is strictly increasing, we may infer that
\[
|x(t)|\le\frac{1}{\varphi_\mu (t)}\quad \forall\,t\in [0,\omega)
\]
and so, {\it a fortiori}, $x(t)\in M:=[-\mu,\mu]$ for all $t\in [0,\omega)$.
Write $K:= P\times M$, $V:= \big[-1,-\half\big]\cup\big[\half,1\big]$ and, as before, define
\[
\chi\colon \R\to\R,~~s\mapsto \min\{vG(\rho,\xi,sv)\mid (\rho,\xi)\in K, ~v\in V\}.
\]
Let $\nu$ be the function
\[
\nu\colon~~\R\to\R_{\ge 0},~~s\mapsto\nu (s):=\max\{\chi(\eta s),0\}.
\]
By Proposition \ref{prop-cont}, $\chi$ is continuous with, by the assumption that $\eta$ is a control direction, $\sup_{s\ge 0}\chi (\eta s)=+\infty$ and so $\nu$ is continuous and surjective.
Fix $\kappa_0 >  \alpha(1/2)$ arbitrarily.  Define an increasing sequence $(\kappa_n)$ in $(\kappa_0,\infty)$ as follows:
choose $\kappa_1 > \kappa_0$ such that $\nu(\kappa_1)=\kappa_0+1$ and set
\[
\kappa_n:=\inf\{\kappa >\kappa_{n-1}\mid \nu (\kappa)=\kappa_0+n\}\quad\forall\, n\ge 2.
\]
\begin{figure}[h]
  \begin{center}
  \begin{tikzpicture}[xscale=.7,yscale=.6]
  \tikzstyle{dot}=[circle, fill, inner sep=1.2pt]
\draw [->, thick,arrows=-Latex] (-1,0)  -- (9.5,0);
\tkzText[right](9.5,0)  {    ~$\kappa$}
\draw [->, thick,arrows=-Latex] (0,-.5)  -- (0,4.5);
\tkzText[above](0.2, 4.5)  { $\nu(\kappa)$};

\draw [dashed] (0,1)  -- (9,1);
\draw [ dashed  ] (0,2)  -- (9,2);
\draw [ dashed] (0,3)  -- (9,3);

\tkzDefPoint(0,1.7){A}

\tkzText[left](0,1)  {      $\kappa_0+1$ \ }
\tkzText[left](0,2)  {      $\kappa_0+2$ \ }
\tkzText[left](0,3)  {      $\kappa_0+3$ \ }
\draw (2.9,1) node[dot] {}; \draw (2.9,0) node[dot] {};
\draw (7,2) node[dot] {}; \draw (8.6,3) node[dot] {}; \draw (7,0) node[dot] {}; \draw (8.6,0) node[dot] {};

 \draw[thick] plot [smooth] coordinates {(0,1.7) (.8,3.5) (1.7,0.0)};
                    \draw[thick] plot [smooth] coordinates {(2.35,0) (3.5,1.7) (4.5,0)};
                     \draw[thick] plot [smooth] coordinates {(5.5,0) (5.9,1) (7,2) (7.45,2.5) (7.9,1.5) (8.6,3) (8.9,3.5)};

 \draw [ dotted, thick ] (2.9,1)  -- (2.9,0);
 \tkzText[below ](3.17,-0.1)  {      $\kappa_1$}

 \draw [ dotted, thick ] (7,2)  -- (7,0);
 \tkzText[below ](7,-0.1)  {      $\kappa_2$}

 \draw [ dotted, thick ] (8.6,3)  -- (8.6,0);
\draw  [ thick ]  (8.6,-0.1) -- (8.6,0.1) ;
 \tkzText[below ](8.6,-0.1)  {      $\kappa_3$}
\end{tikzpicture}
\end{center}
  \vspace{5mm}
  \caption{Schematic construction of the sequence~$(\kappa_n)$.}
  \label{Fig:kappa}
\end{figure}
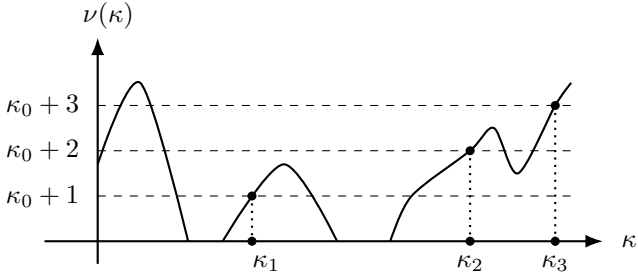
~
\\
Now, we address the question of boundedness of $k(\cdot)$.
\\
Seeking a contradiction, suppose that $k(\cdot)$ is unbounded.
This supposition, together with the observation that  $k(0) =1 <\kappa_0$, ensures that the following is a well-defined sequence~
$(\tau_n)$ in~$(0,\omega)$
\[
\tau_n := \inf\{t\in (0,\omega)\mid k(t)=\kappa_n\},\quad n\in\N.
\]
with the properties
\[
\forall \, n\in \N : \ \tau_{n} < \tau_{n+1}, \quad k(\tau_n ) = \kappa_{n},\quad \nu(k(\tau_n)) =\kappa_0+ n.
\]
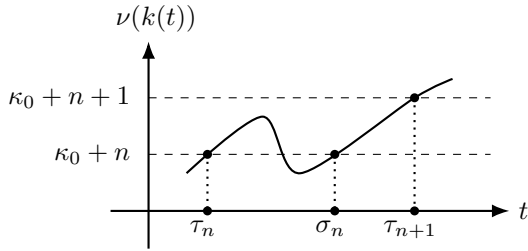
\begin{figure}[h]
  \begin{center}
  \begin{tikzpicture}[xscale=.5,yscale=.5]
  \tikzstyle{dot}=[circle, fill, inner sep=1.2pt]
\draw [->, thick,arrows=-Latex] (-1,0)  -- (9.5,0);
\tkzText[right](9.5,0)  {      $t$}
\draw [->, thick,arrows=-Latex] (0,-1)  -- (0,4.5);
\tkzText[above](0.2, 4.5)  { $\nu(k(t))$};
\tkzText[left](0,1.5)  {      $\kappa_0+n$ \ }
\tkzText[left](0,3)  {      $\kappa_0+n+1$ \ }
\draw (1.55,1.5) node[dot] {}; \draw (1.55,0) node[dot] {};
\draw (4.9,1.5) node[dot] {}; \draw (7,3) node[dot] {}; \draw (4.9,0) node[dot] {}; \draw (7,0) node[dot] {};
\tkzText[below](1.55,0)  {      $\tau_n$ \ }
\tkzText[below](4.9,0)  {      $\sigma_n$ \ }
\tkzText[below](7,0)  {      $\tau_{n+1}$ \ }
  \draw[thick] plot [smooth] coordinates {(1,1)(3,2.5)(4,1)(7,3)(8,3.5)};
 \draw [ dotted, thick ] (1.55,1.5)  -- (1.55,0); \draw [ dotted, thick ] (4.9,1.5)  -- (4.9,0); \draw [ dotted, thick ] (7,3)  -- (7,0);
\draw [dashed] (0,1.5)  -- (9,1.5);
\draw [ dashed] (0,3)  -- (9,3);
\end{tikzpicture}
\end{center}
  \vspace{5mm}
  \caption{Schematic construction of the sequence~$(\sigma_n)$.}
  \label{Fig:sigma}
\end{figure}
Next, define the sequence~$(\sigma_n)$ in $(0,\omega)$ by
\[
\sigma_n:= \sup\{t\in [\tau_n,\tau_{n+1}]\,\mid \nu(k(t))=\kappa_0+n\},\quad n\in\N.
\]
Observe that $\sigma_n < \tau_{n+1}$, $k(\sigma_n) < k(\tau_{n+1})$ and
\[
\kappa_0 + n < \nu (k(t))=\chi (\eta k(t)) \le \nu(k(\tau_{n+1})) = \kappa_0+n+1
\]
for all $t\in (\sigma_n,\tau_{n+1}]$ and all $n\in\N$. Let $n\in\N$ be arbitrary and suppose that, for some $t\in[\sigma_n,\tau_{n+1}]$, $\varphi (t)x(t)\not\in V$.  Then $\varphi(t)|x(t)|<1/2$ and
\[
\alpha (\varphi(t)|x(t)|)=k(t) < \alpha(1/2) < \kappa_0<\kappa_n <\kappa_{n+1}=k(\tau_{n+1})
\]
Therefore, there exists $t^*\in (\sigma_n,\tau_{n+1})$ such that $k(t^*)=\kappa_n$ and we arrive at  the contradiction
\[
\kappa_0+n < \nu(k(t^*))=\nu(\kappa_n)=\kappa_0+n.
\]
Thus, we have shown that
\[
\varphi(t)x(t)\in V\quad \forall\,t\in [\sigma_n,\tau_{n+1}]\quad\forall\,n\in\N.
\]
By symmetry of $V$, we know that $\varphi(t)x(t)\in V$ if, and only if, $-\varphi(t)x(t)\in V$, and so we may infer that
\begin{align*}
-\varphi(t)x(t) &G\big(p(t),x(t),-\eta k(t)\varphi(t)x(t)\big)
\\
&\ \ge \min\{vG(\rho,\xi,\eta k(t)v)\mid (\rho,\xi)\in K,~v\in V\}
\\
&\ =\chi(\eta k(t))
\quad\forall \,t\in [\sigma_n,\tau_{n+1}]\quad \forall\,n\in\N.
\end{align*}
Also (noting that, by bijectivity, $\varphi$ is strictly increasing), we have
\begin{align*}
\varphi(t)\dot\varphi(t)x^2(t) &\le \frac{\dot\varphi(t)}{\varphi(t)}\le c_\varphi \big((1/\varphi (t))+1\big)\\
&\le c_\varphi \big((1/\varphi (\tau_1))+1\big)
=:\gamma   \qquad\forall\, t\in [\tau_1,\omega).
\end{align*}
We may now deduce that
\begin{align*}
\half\tfrac{\mathrm{d}}{\mathrm{d}t}&(\varphi(t)x(t))^2=\varphi(t)\dot\varphi (t)x^2(t)\\
&\qquad\qquad\qquad\quad+\varphi^2(t)x(t)G\big(p(t),x(t),-\eta h(t,x(t)\big)
\\
&\le \gamma -\varphi (t)\big(-\varphi(t)x(t)G\big(p(t),x(t),-\eta k(t)\varphi(t)x(t)\big)\big)
\\
&\le \gamma -\varphi (t)\chi(\eta k(t))
\\
&\le \gamma -\varphi (\tau_1)\nu (k(t))\le \gamma -(\kappa_0 +n)\varphi(\tau_1)
\end{align*}
for all $t\in [\sigma_n,\tau_{n+1}]$ and all $n\in\N$. Choose $n\in\N$ sufficiently large so that $\gamma-(\kappa_0 +n)\varphi(\tau_1) < 0$, in which case
\[
\big(\varphi (\tau_{n+1})x(\tau_{n+1})\big)^2 < \big(\varphi(\sigma_n)x(\sigma_n)\big)^2
\]
and we arrive at the contradiction
\begin{align*}
k(\sigma_n) &< k(\tau_{n+1})=\alpha (\varphi(\tau_{n+1})|x(\tau_{n+1})|)\\
&<\alpha (\varphi(\sigma_n)|x(\sigma_n)|)=k(\sigma_n).
\end{align*}
Therefore, the supposition of unboundedness of $k(\cdot)$ is false.
\end{proof}
The following corollary is immediate, and constitutes the essence of the note. 
\begin{cor}\label{cor}
Let $h$ be as in \eqref{h}.  
If system \eqref{sys-gen} has non-unique control direction, then the control objective is achieved by both the 
feedback $u(t)=-h(t,x(t))$ and its polarity-reversed counterpart  $u(t)=h(t,x(t))$. 
\end{cor}
\section{Example: efficacy of both negative and positive feedback} 
Finally, we proceed to construct an example of a system of form \eqref{sys-gen} with non-unique control direction.  This example, in conjunction with Corollary \ref{cor}, 
establishes the efficacy of a feedback and its polarity-reversed counterpart.    
Let $G$ be given by
\begin{equation}\label{G}
G\colon (\rho,\xi,v)\mapsto f(\rho,\xi) + bg(v),\quad b\ne 0, 
\end{equation}      
with $f\colon\R\to\R$ continuous and 
\begin{equation}\label{g}
g\colon\R\to\R,~~v\mapsto v\sin \big(\ln (1+|v|)\big).
\end{equation}
\begin{prop}\label{techprop}
System \eqref{sys-gen}, with $G$ given by \eqref{G}, has non-unique control direction.  
\end{prop}
\begin{proof}
Let $K\subset\R\times\R$ be compact and $V:= \big[-1,-\half\big]\cup \big[\half,1\big]$.  We will show that, 
for each $\eta\in\{-1,+1\}$, \
\[
\sup_{s\ge 0}\chi (\eta s)=+\infty,
\]
where $\chi\colon s\mapsto \min\{ vf(\rho,\xi)+bvg(sv)\mid (\rho,\xi)\in K,~v\in V\}.$
\\
Define the increasing and unbounded sequence $(s_n)$ by
\[
   s_n :=  \half e^{(n+1)\pi}  - 1 > 0 \quad\forall\,n\in\N,
\]
and so 
\[
\ln (1+s_n) = (n+1)\pi -\ln 2\quad \forall\, n\in\N.
\]
Also, invoking the inequality $e^{\pi/2}> 4$,
\[
1+\half s_n  = \half + \quarter e^{(n+1)\pi} > e^{n\pi}(\quarter e^\pi) > e^{(n+\frac12)\pi}
\]
and so
\[
\ln (1+ \half s_n) > n\pi+\tfrac{\pi}{2}\quad\forall\,n\in\N.
\]
We may now infer
\[
n\pi +\tfrac{\pi}{2} < \ln (1+s_n |v|) \le (n+1)\pi -\ln 2  \quad\forall\, v\in V~~\forall\,n\in\N.
\]
First consider the case of  $n$ even and write $n=2k$.  Since the function $\sin$ is decreasing on the interval $\big[\tfrac{\pi}{2},\pi-\ln 2\big]$, for all $v\in V$ and all $k\in\N$ we have
\begin{align*}
1&=\sin(2k\pi +\tfrac{\pi}{2})  >  \sin(\ln(1+s_{2k}|v|))
\\
&\ge \sin ((2k+1)\pi-\ln 2)=\sin(\pi-\ln 2)=\sin(\ln 2)> 0.
\end{align*}
Writing $\sigma_k := \text{\rm{sgn}}(b)s_k$ for all $k\in\N$, we find
\begin{align*}
bvg(\sigma_{2k}v)&=|b|\text{\rm{sgn}}(b)\,v g(s_{2k}\text{\rm{sgn}}(b)v)
\\
&=s_{2k}|b|v^2\sin (\ln (1+s_{2k}|v|)) \ge\quarter s_{2k}|b|\sin(\ln 2).
\end{align*}
Therefore,  writing 
\[
    c_1:=\min\{vf(\rho,\xi)\mid (\rho,\xi)\in K,~v\in V\},
\]
we have
\begin{align*}
\chi (\sigma_{2k})  &\ge
c_1 + \min\{bvg(\sigma_{2k} v)\mid v\in V\}
\\
&\ge c_1 + \quarter s_{2k}|b|\sin(\ln 2)\quad\forall\,k\in\N.
\end{align*}
For $\eta=\text{\rm{sgn}}(b)$, it follows that $\chi(\eta s_{2k})\to +\infty$ as $k\to\infty$ 
and so $\sup_{s\ge 0} \chi (\eta s) =+\infty$.

Now consider the case of $n$ odd and write $n=2k-1$.  The function $\sin$ is increasing on the interval $\big[-\tfrac{\pi}{2},-\ln 2\big]$ and so
\begin{align*}
 -1&=\sin\big((2k-1)\pi +\tfrac{\pi}{2}\big) < \sin(\ln(1+s_{2k-1}|v|))
 \\
 &\le\sin (2k\pi-\ln 2)= -\sin (\ln 2)\quad\forall\,v\in V~~\forall\,k\in\N.
\end{align*}
Therefore, 
\begin{align*}
bvg(-\sigma_{2k-1}v)&=-|b|s_{2k-1}v^2 \sin (\ln(1+s_{2k-1}|v|))\\
& \ge \quarter |b| s_{2k-1}\sin(\ln2)\quad\forall\,v\in V~~\forall\,k\in\N
\end{align*}
and we have 
\[
\chi (-\sigma_{2k-1}) \ge c_1 + \quarter s_{2k-1}\sin(\ln 2)\quad\forall\,k\in\N.
\]
For $\eta=-\text{\rm{sgn}}(b)$, it follows that $\chi(\eta s_{2k-1})\to\infty$ as $k\to\infty$,
and so $\sup_{s\ge 0} \chi(\eta s) =+\infty$.

In summary, it has been shown that, for each $\eta\in\{-1,+1\}$, $\sup_{s\ge 0}\chi(\eta s)=+\infty$.  
\end{proof}

\section{Discussion}
In succinctly resolving a conjecture by Morse in the 1980s, Nussbaum \cite{3} sparked activity in an area now often referred to as {\em universal adaptive control}.  
In the simplest of contexts, namely, that of scalar linear systems of form \eqref{lin}, a basic question is: ``in the absence of knowledge of the 
control direction $\text{\rm{sgn}}(b)$, does there exist a (reasonably regular) feedback strategy which ensures that the zero state is globally attractive?''.  
As is now well known, this question can be answered in the affirmative: the adaptive strategy 
\[
u(t)=N(k(t))x(t), \quad \dot k(t)=x(t)^2,\quad k(0)=k^0 \ge 0
\]
achieves the control objective if  the function $N$ has the properties: 
\[
\sup_{\kappa >0}\,\frac 1\kappa\int_0^\kappa N(k)\dd k =+\infty\,;\quad \inf_{\kappa >0}\,\frac 1\kappa \int_0^\kappa N (k)\dd k=-\infty.
\]
For example, the function 
$N\colon \kappa\mapsto \kappa \sin\sqrt{\kappa}$ 
suffices.  Loosely speaking, $N$ ensures that the control gain $N\circ k$ can ``probe'' in each control direction and has the potential to take 
arbitrarily large values in magnitude.   Conceptually, the function $g$ in \eqref{G} plays a similar role.  However, the comparison diverges.  
In the adaptive context, the ``probing'' mechanism is a design feature of the controller, whereas, in the context of the present note, the potential 
for ``probing'' is implicit in the system itself: moreover, the adaptive control structure is dynamic, whereas the control structure \eqref{h}  is static 
(albeit time-dependent). 
\section{Conclusion}
A dichotomy of negative and positive feedback has been 
shown to be invalid in general.  This has been done by establishing the existence of a nonlinear example which serves to 
illustrate that ``negative/positive feedback'' need not equate to ``good/bad behaviour''.  In this sense, the claim in this note is purely existential in nature. 
In essence, it simply provides an example in which a feedback and its reversed-polarity counterpart both cause ``good'' behaviour.  The practical utility and quantitative 
behaviour of the feedback employed therein have not been addressed.  The particular example in Section 3 may be open to a criticism of artificiality: 
we stress that its role is only to show that a system of the general form \eqref{sys-gen} can have non-unique control direction, thereby ensuring the existence of a 
stabilizing feedback with the property that its reversed-polarity counterpart is also stabilizing.

\end{document}